\documentclass[a4paper,12pt]{amsart}
\usepackage{amssymb}
\usepackage{ifthen}
\usepackage{graphicx}
\usepackage{float}
\usepackage{caption}
\usepackage{subcaption}
\usepackage{cite}
\usepackage{amsfonts}
\usepackage{amscd}
\usepackage{amsxtra}
\usepackage{color}

\newtheorem{thm}{Theorem}
\newtheorem{cor}{Corollary}
\newtheorem{lem}{Lemma}

\newtheorem{conj}{Conjecture}
\newtheorem{prob}{Problem}

\theoremstyle{definition}
\newtheorem{defn}{Definition}[section]
\newtheorem{example}{Example}

\newenvironment{pf}[1][]{%
 \vskip 1mm
 \noindent
 \ifthenelse{\equal{#1}{}}%
  {{\slshape Proof. }}%
  {{\slshape #1.} }%
 }%
{\qed\bigskip}

\newcounter{alphabet}
\newcounter{tmp}
\newenvironment{Thm}[1][]{\refstepcounter{alphabet}%
\bigskip%
\noindent%
{\bf Theorem \Alph{alphabet}}%
\ifthenelse{\equal{#1}{}}{}{ (#1)}%
{\bf .} \itshape}{\vskip 8pt}

\makeatletter
\newcommand{\Ref}[1]{\@ifundefined{r@#1}{}{\setcounter{tmp}{\ref{#1}}\Alph{tmp}}}
\makeatother

\newcommand{\IC}{{\mathbb C}}
\newcommand{\ID}{{\mathbb D}}

\def\be{\begin{equation}}
\def\ee{\end{equation}}

\newcommand{\bee}{\begin{enumerate}}
\newcommand{\eee}{\end{enumerate}}

\newcommand{\blem}{\begin{lem}}
\newcommand{\elem}{\end{lem}}
\newcommand{\bthm}{\begin{thm}}
\newcommand{\ethm}{\end{thm}}
\newcommand{\bcor}{\begin{cor}}
\newcommand{\ecor}{\end{cor}}
\newcommand{\beg}{\begin{example}}
\newcommand{\eeg}{\end{example}}
\newcommand{\begs}{\begin{examples}}
\newcommand{\eegs}{\end{examples}}
\newcommand{\bdefe}{\begin{defn}}
\newcommand{\edefe}{\end{defn}}
\newcommand{\bprob}{\begin{prob}}
\newcommand{\eprob}{\end{prob}}
\newcommand{\bques}{\begin{ques}}
\newcommand{\eques}{\end{ques}}
\newcommand{\bei}{\begin{itemize}}
\newcommand{\eei}{\end{itemize}}
\newcommand{\bcon}{\begin{conj}}
\newcommand{\econ}{\end{conj}}
\newcommand{\bcons}{\begin{conjs}}
\newcommand{\econs}{\end{conjs}}
\newcommand{\bprop}{\begin{propo}}
\newcommand{\eprop}{\end{propo}}
\newcommand{\br}{\begin{rem}}
\newcommand{\er}{\end{rem}}
\newcommand{\brs}{\begin{rems}}
\newcommand{\ers}{\end{rems}}
\newcommand{\bo}{\begin{obser}}
\newcommand{\eo}{\end{obser}}
\newcommand{\bos}{\begin{obsers}}
\newcommand{\eos}{\end{obsers}}
\newcommand{\bpf}{\begin{pf}}
\newcommand{\epf}{\end{pf}}
\newcommand{\ba}{\begin{array}}
\newcommand{\ea}{\end{array}}
\newcommand{\beq}{\begin{eqnarray}}
\newcommand{\beqq}{\begin{eqnarray*}}
\newcommand{\eeq}{\end{eqnarray}}
\newcommand{\eeqq}{\end{eqnarray*}}

\newcommand{\ds}{\displaystyle}

\newcounter{minutes}
\divide\time by 60
\newcounter{hours}
\multiply\time by 60 \addtocounter{minutes}{-\time}

\begin{document}
\bibliographystyle{amsplain}
\title[On a powered Bohr inequality]{On a powered Bohr inequality}

\thanks{
File:~\jobname .tex,
          printed: \number\day-\number\month-\number\year,
          \thehours.\ifnum\theminutes<10{0}\fi\theminutes}


\author{Ilgiz R Kayumov, and Saminathan Ponnusamy }

\address{I. R Kayumov, Kazan Federal University, Kremlevskaya 18, 420 008 Kazan, Russia
}
\email{ikayumov@kpfu.ru }

\address{S. Ponnusamy, Department of Mathematics,
Indian Institute of Technology Madras, Chennai--600 036, India.}
\email{samy@isichennai.res.in, samy@iitm.ac.in}


\begin{abstract}
The object of this paper is to study the powered Bohr radius $\rho_p$, $p \in (1,2)$, of analytic functions
$f(z)=\sum_{k=0}^{\infty} a_kz^k$ and such that $|f(z)|<1$ defined on the unit disk $|z|<1$. More precisely, if
$M_p^f (r)=\sum_{k=0}^\infty |a_k|^p r^k$, then we show that $M_p^f (r)\leq 1$ for  $r \leq r_p$ where $r_\rho$ is the powered Bohr radius
for conformal automorphisms of the unit disk.
This answers the open problem posed by Djakov and Ramanujan  in 2000. A couple of other consequences of our
approach is also stated, including an asymptotically sharp form of one of the
results of Djakov and Ramanujan. In addition, we consider a similar problem for sense-preserving harmonic mappings in $|z|<1$.
Finally, we conclude by stating the Bohr radius for the class of Bieberbach-Eilenberg functions.
\end{abstract}

\thanks{
}

\maketitle
\pagestyle{myheadings}
\markboth{I. R. Kayumov and S. Ponnusamy}{On a powered Bohr inequality}

\section{Preliminaries and Main Results}
Let $\mathcal B$ denote the class of analytic functions $f$ defined on the unit disk $\ID :=\{z\in\IC:\, |z|<1\}$, with the power
series expansion $f(z)=\sum_{k=0}^{\infty} a_kz^k$ and such that $|f(z)|<1$ for $z\in\ID$. Then the classical Bohr's inequality states that
there is a constant $\rho$ such that
$$M^f(r):=\sum_{k=0}^{\infty}|a_k|r^k \leq 1 
~\mbox{ for all $r=|z|\leq \rho$}
$$
and the value $\rho=1/3$ is optimal. The number $\rho=1/3$, known as Bohr's radius, was originally obtained in 1914 by H.~Bohr \cite{Bohr-14} with
$\rho=1/6$, but subsequently later, Wiener, Riesz and Schur, independently established
the sharp inequality for $r=|z|\leq 1/3$. This little article of Bohr generates intensive research activities even after a century of its appearance.
We refer to the recent survey article on this topic \cite{AAPon1} and the references therein.
Multidimensional generalizations of this result were obtained by Boas and Khavinson \cite{BoasKhavin-97-4} by establishing upper and lower bounds for the
Bohr radius of the unit polydisk $\ID^n$.
Aizenberg \cite{Aizen-00-1,Aizen-05-3} extended the concept of Bohr radius in several different directions for further studies in this topic.
In 2000, Djakov and Ramanujan \cite{DjaRaman-2000} investigated
the same phenomenon from different point of view. For  $f\in {\mathcal B}$ and a fixed $p>0$, we consider the powered Bohr sum $M_p^f (r)$ defined by
%
$$
M_p^f (r)=\sum_{k=0}^\infty |a_k|^p r^k.
$$
Observe that for $p=1$, $M_p^f (r)$ reduces to the classical Bohr sum defined as above by $M^f(r)$. The best possible constant $\rho_p$ for which
$$M_p^f (r)\leq 1 ~\mbox{ for all $r\leq \rho_p$}
$$
is called the (powered) Bohr radius for the family ${\mathcal B}$.

We now introduce
$$
M_p(r):=\sup_{f \in {\mathcal B}} M_p^f(r)
$$
and
$$ r_p:=\sup\left \{r:\, a^p +\frac{r(1-a^2)^p}{1-r a^p} \leq 1, \quad 0 \leq a <1\right \}=\inf_{a \in [0,1)} \frac{1-a^p}{a^p(1-a^{p})+(1-a^2)^p}.
$$
Let us first proceed to recall the following results.

\begin{Thm}{\rm (\cite[Theorem 3]{DjaRaman-2000})}
For each  $p \in (1,2)$ and $f(z)=\sum_{k=0}^{\infty} a_kz^k$ belongs to $\mathcal B$, we have
$M_p^f (r)\leq 1$ for  $r \leq T_p$, where
$$ m_p\leq T_p\leq r_p.
$$
Here $r_p$ is as above and
$$m_p:=\frac{p}{\left (2^{1/(2-p)}+p^{1/(2-p)}\right )^{2-p}}.
$$
\end{Thm}

\begin{Thm}{\rm (\cite[Theorem 2]{DjaRaman-2000})}
For each  $p \in (0,2)$
$$
M_p(r) \asymp \left(\frac{1}{1-r}\right)^{1-p/2}.
$$
\end{Thm}

Our first aim is to investigate the problem posed by Djakov and Ramanujan \cite{DjaRaman-2000} about the Bohr radius for $M_p^f(r)$.
Their question is the following.

\bprob {\rm \cite[Question 1, p.~71]{DjaRaman-2000}}
What is the exact value of the (powered) Bohr radius $\rho_p$, $p\in (1,2)$? Is it true that $\rho_p=r_p$?
\eprob

Using the method of proofs of our recent approach from \cite{KayPon1,KayPon2}, we solve this problem affirmatively in the following form.

\begin{thm}\label{KP9-thm1}
If $f(z)=\sum_{k=0}^{\infty} a_kz^k$ belongs to $\mathcal B$ and $0<p \leq 2$, then
$$M_p(r)=\max_{a \in [0,1]}\left[ a^p +\frac{r(1-a^2)^p}{1-r a^p}\right], \quad 0 \leq r \leq 2^{p/2-1},
$$ and
$$
M_p(r) < \left(\frac{1}{1-r^{2/(2-p)}}\right)^{1-p/2}, \quad   2^{p/2-1} < r < 1.
$$
\end{thm}

Proofs of Theorem \ref{KP9-thm1} and a couple of its corollaries will be given in Section \ref{KayPon9-sec2}.

Let us remark that $M_p(r)=1$ for $p \ge 2$ and $r \leq 1$. So, the interesting case is to consider the problem only for
$p \in (1,2)$.

One may ask about the second inequality of Theorem \ref{KP9-thm1}: how close it to be sharp?
To get an answer to this question we will use a Bombieri-Bourgain estimate \cite{BombBour-2004} which reads as follows: for a given
$\varepsilon >0$,
there exists a positive constant $C(\varepsilon)>0$, such that
$$ M_1(\rho) \ge \frac{1}{\sqrt{1-\rho^2}}-C(\varepsilon)\left ( \log \frac{1}{1-\rho}\right )^{(3/2)+\varepsilon}, \quad \rho \ge 1/\sqrt{2}.
$$
The H\"{o}lder inequality implies that
\beqq
M_1(f,r^{1/(2-p)}) &=&\sum_{k=0}^\infty |a_k| r^{k/p}r^{(2k(p-1))/(p(2-p))}\\
&\leq & \left(\sum_{k=0}^\infty |a_k|^p r^k\right)^{1/p}\left(\sum_{k=0}^\infty r^{2k/(2-p)}\right)^{1-1/p}\\
&=& \left (M_p^f(r)\right )^{1/p}\frac{1}{(1-r^{2/(2-p)})^{(p-1)/p}}
\eeqq
so that
$$
M_p^f(r) \ge  \left(\frac{1}{\sqrt{1-r^{2/(2-p)}}}-C(\varepsilon)\left (\log \frac{1}{1-r^{1/(2-p)}}\right )^{3/2+\varepsilon} \right)^p(1-r^{2/(2-p)})^{p-1}, \quad   2^{p/2-1} < r < 1,
$$ or equivalently
$$
M_p^f(r) \ge \left(\frac{1}{1-r^{2/(2-p)}}\right)^{1-p/2}-C_1(\varepsilon)(1-r^{2/(2-p)})^{(p-1)/2}\left (\log \frac{1}{1-r^{1/(2-p)}}\right )^{3/2+\varepsilon}.
$$
This estimate together with the second estimate of Theorem \ref{KP9-thm1} implies that
$$M_p(r) -  \left(\frac{1}{1-r^{2/(2-p)}}\right)^{1-p/2} \to 0 ~\mbox{ as }~ r \to 1^{-}
$$
for $1<p<2$ while we do not know whether this fact is true for $p=1$. Also the last estimate can be considered as an
asymptotically sharp form of Theorem B in the case $p>1$.

\begin{cor}\label{KP9-cor1}
Let $p \in (1,2)$. Then $M_p(r) = 1$ for $r\leq r_p$.
\end{cor}

In \cite[Corollary 2.8]{PaulPopeSingh-02-10}, Paulsen et al. showed that if $f\in {\mathcal B}$, then for $r\in [0,1)$,
\be\label{KP9-eq5}
M_1^f (r)\leq m(r) =\inf\{M(r), 1/\sqrt{1-r^{2}}\}
\ee
where
$$M(r)=  \sup\left \{ t +( 1-t^2)\frac{r}{1-r}:\, 0\leq t\leq 1\right \}
=\left \{\begin{array}{cl}
1& \mbox{ for }~ 0\leq r\leq 1/3 \\
\ds \frac{4r^{2}+(1-r)^2}{4r(1-r)} & \mbox{ for }~ 1/3<r<1.
\end{array}\right .
$$
In 2002,   Paulsen et al. \cite{PaulPopeSingh-02-10} raised a question whether the inequality \eqref{KP9-eq5} is sharp for any
$r$ with $1/3<r<1$. However, in 1962 this has been answered by Bombieri \cite{BombBour-196} who determined
the exact value of this constant for $r$ in the range $1/3\leq  r\leq 1/\sqrt{2}.$ This constant is
$$m(r ) =\frac{3-\sqrt{8(1-r^2)}}{r}.
$$
Further results on this and related topics can be found in \cite{DjaRaman-2000,PaulPopeSingh-02-10}. On the other hand, it is worth mentioning that
the answer to the above question is indeed a consequence of Theorem \ref{KP9-thm1} and so, we state it as a corollary.

\begin{cor}\label{KP9-cor2}
We have the following sharp estimate:
$$M_1(r) = \frac{1}{r}(3-\sqrt{8(1-r^2)}) \mbox{ for } r \in \left[\frac{1}{3}, \frac{1}{\sqrt{2}} \right].
$$
\end{cor}

Finally, we recall the following corollary which was proved in \cite{KayPon2} and so we omit the proof.
\begin{cor}\label{KP2-th3}
Let $p \in \mathbb{N}$ and $0\leq m \leq p$, $f(z)=\sum_{k=0}^{\infty} a_{pk+m}z^{pk+m}$ be analytic in $\ID$
and $|f(z)| < 1$ in $\ID$. Then
$$\sum_{k=0}^{\infty} |a_{pk+m}|r^{pk+m} \leq 1 \mbox{ for } r \leq r_{p,m},
$$
where $ r_{p,m}$ is the maximal positive root of the equation
$$
-6 r^{p-m} + r^{2(p-m)} + 8 r^{2p} +1= 0.
$$
 The extremal function has the form $z^m(z^p-a)/(1-az^p)$, where
$$a=\left (1-\frac{\sqrt{1-{r_{p,m}}^{2p}}}{\sqrt{2}}\right )\frac{1}{{r_{p,m}}^p }.
$$
\end{cor}

Our next result concerns sense-preserving harmonic mappings defined on the unit disk $\ID$. Recall that the
family ${\mathcal H}$ of complex-valued harmonic functions $f=h+\overline{g}$ defined on ${\mathbb D}$ and its univalent subfamilies
are investigated in details. Here $h$ and $g$ are analytic on $\ID$ with the form
$$
h(z)=\sum_{k=0}^{\infty}a_kz^k~\mbox{ and }~g(z)=\sum _{k=1}^{\infty}b_kz^k
$$
so that the Jacobian of $f$ is given by $J_f= |f_z|^2-|f_{\overline{z}}|^2=|h'|^2-|g'|^2$. We say that
the harmonic mapping $f$ is sense-preserving if $J_f(z)>0$ in $\ID$.
We call $\omega (z)=g'(z)/h'(z)$ the complex dilatation of
$f=h+\overline{g}$. Lewy's theorem implies that every harmonic
function $f$ on $\ID$ is locally one-to-one and sense-preserving on
$\ID$ if and only if $|\omega (z)|<1$ for $z\in \ID$. See \cite{Clunie-Small-84,D}
for detailed discussion on the class of univalent harmonic mappings and its geometric subclasses.

\begin{thm}\label{KayPon9-thm2}
Suppose that $f(z) = h(z)+\overline{g(z)}=\sum_{k=0}^\infty a_k z^k+\overline{\sum_{k=1}^\infty b_k z^k}$ is a
harmonic mapping of the disk $\ID$, where $h$ is a bounded function in $\ID$ and $|g'(z)|\leq |h'(z)|$ for $z\in\ID$
(the later condition obviously holds if $f$ is sense-preserving). If $p \in [0,2]$ then the following sharp inequality holds
$$
|a_0|^p+\sum_{k=1}^\infty (|a_k|^p+ |b_k|^p) r^k \leq  ||h||_{\infty}\max_{a \in [0,1]}\left\{ a^p +\frac{2r(1-a^2)^p}{1-r a^p}\right\}
$$
for  $\ds r \leq (2^{1/(p-2)}+1)^{p/2-1}$.
In the case $p>2$ we have
$$|a_0|^p+\sum_{k=1}^\infty (|a_k|^p+ |b_k|^p) r^k \leq  ||h||_{\infty} \max\{1,2r\}.
$$
\end{thm}

\begin{cor}\label{KayPon9-cor4}
Suppose that $f(z) = h(z)+\overline{g(z)}=\sum_{k=0}^\infty a_k z^k+\overline{\sum_{k=1}^\infty b_k z^k}$ is a
sense-preserving harmonic mapping of the disk $\ID$, where $h$ is a bounded function in $\ID$. Then the following sharp inequalities holds:
$$|a_0|+\sum_{k=1}^\infty (|a_k|+ |b_k|) r^k \leq  \frac{||h||_{\infty}}{r}(5-2\sqrt{6}\sqrt{1-r^2})~\mbox{ for  }~ \frac{1}{5} \leq r \leq \sqrt{\frac{2}{3}},
$$
and
$$
|a_0|+\sum_{k=1}^\infty (|a_k|+ |b_k|) r^k \leq  ||h||_{\infty}~\mbox{ for  }~ r \leq \frac{1}{5}.
$$
\end{cor}

Proofs of Theorem \ref{KayPon9-thm2} and Corollary \ref{KayPon9-cor4} will be given in Section \ref{KayPon9-sec2}. In Section \ref{KayPon9-sec3},
we discuss Bohr radius for the class of Bieberbach-Eilenberg functions.

\section{Proofs of Theorems \ref{KP9-thm1} and \ref{KayPon9-thm2} and their corollaries}\label{KayPon9-sec2}

The proofs of the theorems rely on a couple of lemmas established by the present authors in \cite{KayPon1} (see also \cite{KayPon2}).

\begin{lem}\label{KP2-lem2}
{\rm \cite{KayPon1}}
 Let $|a|<1$ and $0 < R \leq 1$. If $g(z)=\sum_{k=0}^{\infty} b_kz^k$ belongs to $\mathcal B$, then
the following sharp inequality holds:
$$
\sum_{k=1}^\infty |b_k|^2R^{k} \leq R\frac{(1-|b_0|^2)^2}{1-|b_0|^2R}.
$$
\end{lem}

\begin{lem}\label{lem22} For all $p \in (0,2)$, we have $ r_p < \left(1/2\right)^{1-p/2}. $
\end{lem}
\bpf Let $r=r_p$ and set $a=\left(1/2\right)^{1-p/2}$. Then we conclude that
$$
a^p+r\frac{(1-a^2)^p}{1-ra^p} =2\left(\frac{1}{2}\right)^{p/2}>1
$$ which contradicts to the definition of $r_p$.
\epf


\bpf[Proof of Theorem \ref{KP9-thm1}]
Let $|a_0|=a>0$ and $r \leq 2^{p/2-1}$. At first we suppose that $a > r^{1/(2-p)}$. In this case we have
\beqq
M_p^f (r) & = & a^p+\sum_{k=1}^\infty \rho^k|a_k|^p\left(\frac{r}{\rho}\right)^k\\
& \leq & a^p+\left(\sum_{k=1}^\infty \left(\rho^k|a_k|^p\right)^{2/p}\right)^{p/2}\left(\sum_{k=1}^\infty \left(\frac{r}{\rho}\right)^{2k/(2-p)}\right)^{1-p/2}\\
&=&a^p+\left(\sum_{k=1}^\infty (\rho ^{2/p})^{k}|a_k|^2\right)^{p/2}\left(\sum_{k=1}^\infty \left (\left(\frac{r}{\rho}\right)^{2/(2-p)}\right )^k\right)^{(2-p)/2}\\
&\leq &a^p+ \left ( \frac{\rho^{2/p} (1-a^2)^2}{1-a^2\rho^{2/p}}\right )^{p/2}\left(\frac{ (r/\rho)^{2/(2-p)}}{1-(r/\rho)^{2/(2-p)}}\right)^{(2-p)/2}
~\mbox{ (by Lemma \ref{KP2-lem2}),}\\
&=&a^p+r \left(\frac{(1-a^2)^2}{1-a^2\rho^{2/p}}\right)^{p/2}\left(\frac{1}{1-(r/\rho)^{2/(2-p)}}\right)^{(2-p)/2}.
\eeqq
Setting $\rho=r^{p/2} a^{(p-2)p/2}$ we obtain the inequality
$$
M_p^f (r) \leq a^p+r\frac{(1-a^2)^p}{1-ra^p},
$$
which proves the theorem in the case $a > r^{1/(2-p)}$.

In the case $a \leq r^{1/(2-p)}$, we set $\rho = 1$ and obtain
$$M_p^f (r)=\sum_{k=0}^\infty |a_k|^p r^k \leq a^p +r\frac{(1-a^2)^{p/2}}{(1-r^{2/(2-p)})^{1-p/2}}.
$$
Let us remark that the inequality $M_p^f (r) \leq 1$ is valid in the cases $a=0$ and $a=r^{1/(2-p)}$. This fact can be established
as a limiting case of the previous case.
Finally, we let $t=a^2$. We have then to maximize the expression
$$ A(t)=t^{p/2} +r\frac{(1-t)^{p/2}}{(1-r^{2/(2-p)})^{1-p/2}}, \quad t \leq r^{2/(2-p)}.
$$
Using differentiation we obtain the stationary point
$$t=1-r^{2/(2-p)}
$$
which must satisfy under the restriction $t\leq r^{2/(2-p)}$ which is impossible because $r \leq 2^{p/2-1}$.

However, in the case $r > 2^{p/2-1}$ the critical point $t$ is admissible so that
$$A(t)=t^{p/2} +r\frac{(1-t)^{p/2}}{(1-r^{2/(2-p)})^{1-p/2}}=\left(\frac{1}{1-r^{2/(2-p)}}\right)^{1-p/2}.
$$
This observation shows that
$$
M_p^f(r) \leq \left(\frac{1}{1-r^{2/(2-p)}}\right)^{1-p/2}, \quad   2^{p/2-1} < r < 1.
$$ Now let us show that this inequality cannot be sharp. To do this we will use the method presented
by Bombieri and Bourgain \cite{BombBour-2004}.

Suppose that the estimate sharp in this case. Then by analyzing H\"{o}lder's inequality we immediately conclude that
$$|a_k|=\sqrt{1-r^{2/(2-p)}}\, r^{k/(2-p)}, \quad k \ge 0.
$$
Also it is easy to show that the extremal function must be a Blashke's product with a finite degree $d \ge 1$. Computing the
area, one obtains that 
$$
\pi d = {\rm Area}\,f(\ID)\, =\pi\sum_{k=1}^\infty k |a_k|^2= \pi \frac{\lambda^2}{1-\lambda^2}, \quad \lambda=r^{1/(2-p)}.
$$ 
From here we easily deduce that $d=\lambda^2/(1-\lambda^2)$ and thus, $\lambda =\sqrt{d/(d+1))}$, which gives
\begin{equation}\label{extra1}
\sqrt{\frac{d}{d+1}} =r^{1/(2-p)}, ~\mbox{ i.e. }~ r=\left(\frac{d}{d+1}\right)^{1-(p/2)}.
\end{equation}
Therefore our inequality could be sharp for these values only. Now let us show that this is possible for $d=1$ only.
Using the same reasoning as in \cite{BombBour-2004} (in fact we apply their considerations in which $r$ is replaced by $r^{1/(2-p)}$)
we arrive at the identity
$$
\sqrt{1-r^{2/(2-p)}}=r^{d/(2-p)}
$$ which together with (\ref{extra1}) implies that
$$
\sqrt{1-\frac{d}{d+1}}= \left(\frac{d}{d+1}\right)^{d/2}.
$$
From here we easily deduce that $d=1$ and this completes the proof of Theorem \ref{KP9-thm1}.
\epf

\bpf[Proof of Corollary \ref{KP9-cor1}] Easily follows from Theorem  \ref{KP9-thm1} and Lemma \ref{lem22}.
\epf

\bpf[Proof of Corollary \ref{KP9-cor2}]
Theorem \ref{KP9-thm1} for $p=1$ gives that
$$ M_p(r)=\max_{a \in [0,1]}\left\{ a +\frac{r(1-a^2)}{1-r a}\right \}.
$$ 
By using differentiation it is easy to show that in the case $1/3 \leq r \leq 1/\sqrt{2}$ the 
maximum of the last expression is achieved at the point
$$
a=\left (1-\frac{\sqrt{1-{r}^{2}}}{\sqrt{2}}\right )\frac{1}{{r}}
$$
and consequently, we obtain that
$$ M_1(r)=\frac{1}{r}(3-2\sqrt{2}\sqrt{1-r^2}).
$$
The proof is complete.
\epf

\bpf[Proof of Theorem \ref{KayPon9-thm2}]
Without lost of generality we may assume that $||h||_\infty=1$. As in \cite{KayPon3}, the condition 
$|g'(z)|\leq |h'(z)|$ gives that for each $r\in [0,1)$,
\be\label{KP9-eq6}
\sum_{k=1}^\infty |b_k|^2r^k \leq \sum_{k=1}^\infty |a_k|^2r^k.
\ee
Let $|a_0|=a>0$. Then, by using the same method as in the previous theorem in the case $a > r^{1/(2-p)}$, we obtain
$$
|a_0|^p+\sum_{k=1}^\infty (|a_k|^p+ |b_k|^p) r^k \leq a^p+2r\frac{(1-a^2)^p}{1-ra^p}.
$$

In the case $a \leq r^{1/(2-p)}$, we let $\rho = 1$ and obtain
$$\sum_{k=0}^\infty |a_k|^p r^k \leq a^p +2r\frac{(1-a^2)^{p/2}}{(1-r^{2/(2-p)})^{1-p/2}}.
$$

We set $t=a^2$. We have to maximize the expression
$$ B(t)=t^{p/2} +2r\frac{(1-t)^{p/2}}{(1-r^{2/(2-p)})^{1-p/2}}, \quad t \leq r^{2/(2-p)}.
$$
Using differentiation we see that the function $B(t)$ is increasing on the interval
$$0 \leq t \leq \frac{1-r^{2/(2-p)}}{1+(2r)^{2/(2-p)}-r^{2/(2-p)}}.
$$
The upper bound of this interval is greater than or equal to $2^{p/2-1}$ in the case $r \leq (2^{1/(p-2)}+1)^{p/2-1}$.
It means that the function $B(t)$ has maximum at the point $t=r^{2/(2-p)}$ which corresponds to the case $a=r^{1/(2-p)}$
so that we can apply our previous case. This completes the proof Theorem \ref{KayPon9-thm2}.
\epf

Let $p=1$ and then we apply the previous theorem. As a result, we obtain the inequality
$$
|a_0|+\sum_{k=1}^\infty (|a_k|+ |b_k|) r^k \leq  \max_{a \in [0,1]}\left\{ a +\frac{2r(1-a^2)}{1-r a}\right\}~\mbox{ for  }~ r \leq \sqrt{2/3}.
$$
Straightforward calculations confirm the proof of Corollary \ref{KayPon9-cor4}.
In Section \ref{KayPon9-sec3}, we present the Bohr radius for the class of Bieberbach-Eilenberg functions.

\section{Concluding remarks} \label{KayPon9-sec3}
Let $\mathcal{BE}$ denote the class of all functions $f(z)=\sum_{k=1}^{\infty} a_kz^k$ analytic in $\ID$ such that $f(z_1)f(z_2)\neq 1$ for all pairs
of points $z_1,z_2$ in $\ID$. Each $f\in\mathcal{BE}$ is called a Bieberbach-Eilenberg function. Clearly, $\mathcal{BE}$  contains the class $\mathcal{B}_0$,
where $\mathcal{B}_0=\{f\in \mathcal{B}:\, f(0)=0\}.$ In 1970, Aharonov \cite{Ahar-70} and Nehari \cite{Nehar-70} independently showed that
\be\label{KP9-eq7}
\sum_{k=1}^\infty |a_k|^2 \leq 1 ~\mbox{ and }~ |f(z)| \leq \frac{|z|}{\sqrt{1-|z|^2}}
\ee
hold for every $f\in\mathcal{BE}$. Equality holds only for the functions
$$ f(z)=\frac{\eta z}{R\pm (\sqrt{R^2-1})i\eta z}, \quad R>1,~|\eta|=1.
$$
Since $\mathcal{B}_0\subset \mathcal{BE}$, it is natural to ask for the Bohr radius for the family $\mathcal{BE}$.
Indeed, we see blow that the Bohr radius for $\mathcal{BE}$ and the class $\mathcal{B}_0$ remains the same.

\bthm\label{KP9-thm3}
Assume that $f(z)=\sum_{k=1}^{\infty} a_kz^k$ belongs to $\mathcal{BE}$. Then
$$
\sum_{k=1}^\infty |a_k|r^k \leq 1 ~\mbox{ for $|z|=r \leq 1/\sqrt{2}$}.
$$
The number $1/\sqrt{2}$ is sharp.
\ethm
\bpf
Because $f\in\mathcal{BE}$ satisfies the coefficient inequality \eqref{KP9-eq7}, it follows that
$$\sum_{k=1}^\infty |a_k|r^{k} \leq \sqrt{\sum_{k=1}^\infty |a_k|^2 }{\sqrt{\sum_{k=1}^\infty r^{2k}}} \leq  \frac{r}{\sqrt{1-r^2}}
$$
which is less than or equal to $1$ if $0\leq r\leq 1/\sqrt{2}$. The number $1/\sqrt{2}$ is sharp as the function
$f(z)=z(a-z)/(1-az)$ shows, where $a=1/\sqrt{2}$. The proof is complete.
\epf

\begin{thm}\label{KayPon9-thm3}
Suppose that $f(z) = h(z)+\overline{g(z)}=\sum_{k=1}^\infty a_k z^k+\overline{\sum_{k=1}^\infty b_k z^k}$ is a
harmonic mapping of the disk $\ID$, where $h\in\mathcal{BE}$ and $|g'(z)|\leq |h'(z)|$ for $z\in\ID$.
Then for any $p \ge 1$ and $r<1$, the following inequality holds:
\beqq
\sum_{k=1}^\infty(|a_k|^p +|b_k|^p)^{1/p}r^{k}  &\leq &
 \max \{2^{(1/p)-1/2},1\}\frac{\sqrt{2}r}{\sqrt{1-r^2}}.
\eeqq
\end{thm}
\bpf
By hypothesis, \eqref{KP9-eq6} holds and thus,  letting $r$ approach $1$, we get
$$\sum_{k=1}^\infty |b_k|^2 \leq \sum_{k=1}^\infty |a_k|^2 \leq 1.
$$
Consequently, we obtain
\beqq
\sum_{k=1}^\infty (|a_k|^p +|b_k|^p)^{1/p}r^{k} &\leq &\sqrt{\sum_{k=1}^\infty(|a_k|^p +|b_k|^p)^{2/p}}{\sqrt{\sum_{k=1}^\infty r^{2k}}}\\
&\leq& \sqrt{\max\{ 2^{(2/p)-1}, 1\}\sum_{k=1}^\infty(|a_k|^2 +|b_k|^2)}\, \frac{r}{\sqrt{1-r^2}}\\
& \leq & \max\{2^{(1/p)-1/2},1\}\frac{\sqrt{2}\, r}{\sqrt{1-r^2}}
\eeqq
and the proof is complete.
\epf

Theorem \ref{KayPon9-thm3} for $p=1$ shows that for $r\leq 1/\sqrt{5}$,
$$\sum_{k=1}^\infty (|a_k| +|b_k|) r^k \leq 1.
$$
Similarly, for $p =2$, we see that for $r\leq 1/\sqrt{3}$,
$$\sum_{k=1}^\infty (|a_k|^2 +|b_k|^2)^{1/2}r^{k} \leq 1.
$$

\subsection*{Acknowledgements}
The research of the first author was supported by Russian foundation for basic research, Proj. 17-01-00282.
The work of the second author is supported  in part by Mathematical Research Impact Centric Support
(MATRICS) grant, File No.: MTR/2017/000367, by the Science and Engineering Research Board (SERB),
Department of Science and Technology (DST), Government of India.

\end{document}